\documentclass[11pt]{article}%
\usepackage{amsmath}
\usepackage{amsfonts}
\usepackage{verbatim}
\usepackage{amssymb}
\usepackage{graphicx}
\usepackage{a4wide}
\usepackage{hyperref}%
\providecommand{\U}[1]{\protect\rule{.1in}{.1in}}
\newtheorem{theorem}{Theorem}

\newtheorem{corollary}[theorem]{Corollary}

\newtheorem{lemma}[theorem]{Lemma}

\newtheorem{proposition}[theorem]{Proposition}

\begin{document}

\title{Characterizations of objective sets and objective functions}
\author{C. Z\u{a}linescu\\University Alexandru Ioan Cuza Ia\c{s}i, Faculty of Mathematics \\700506 Ia\c{s}i, Romania (email: zalinesc@uaic.ro).}
\date{}
\maketitle

\textbf{Abstract.} In this short note we give characterizations for objective
sets and objective functions as defined in D.Y. Gao and C. Wu's paper arxiv:1104.2970v2.

\section{Objective sets and objective functions}

In \cite{GaoWu:11v2} one can read the following:

\textquotedblleft Mathematical definitions of the objective set and objective
function are given in the book [9] (Chapter 6, page 288). Let

$\mathcal{Q}=\{Q\in\mathbb{R}^{m\times m}\mid Q^{T}=Q^{-1},\ \det Q=1\}$

\noindent be a proper orthogonal rotation group.

\textbf{Definition 1} \textbf{(Objectivity and Isotropy)} A subset
$\mathcal{Y}_{a}\subset\mathbb{R}^{m}$ is said to be \emph{objective} if
$Qy\in\mathcal{Y}_{a}$ $\forall y\in\mathcal{Y}_{a}$ and $\forall
Q\in\mathcal{Q}$. A real-valued function $T:\mathcal{Y}_{a}\rightarrow
\mathbb{R}$ is said to be \emph{objective} if its domain is objective and

$T(Qy)=T(y)$ $\forall y\in\mathcal{Y}_{a}$ and $\forall Q\in\mathcal{Q}$.
$(2)$

A subset $\mathcal{Y}_{a}\subset\mathbb{R}^{m}$ is said to be \emph{isotropic}
if $yQ^{T}\in\mathcal{Y}_{a}$ $\forall y\in\mathcal{Y}_{a}$ and $\forall
Q\in\mathcal{Q}$. A real-valued function $T:\mathcal{Y}_{a}\rightarrow
\mathbb{R}$ is said to be \emph{isotropic} if its domain is isotropic and

$T(yQ^{T})=T(y)$ $\forall y\in\mathcal{Y}_{a}$ and $\forall Q\in Q$.
$(3)$\textquotedblright.

\smallskip

The reference ``[9]'' above is our reference \cite{Gao:00}. Because $Q$ is an
$m\times m$ (orthogonal) matrix, in order to calculate $Qy$ for $y\in
\mathbb{R}^{m}$ one must consider $y$ as a column vector. Of course, $Q^{T}$
is also an $m\times m$-matrix. Because we don't know what is meant by $yQ^{T}$
in such a case, in the sequel we deal only with objective sets and objective
functions in the sense mentioned above.

\smallskip

Let us also see how these notions are defined in Definition 6.1.2 in
\cite{Gao:00}:

\textquotedblleft\textbf{Definition 6.1.2 (Objectivity and Isotropy)}

(D1) \emph{Objective Set and Objective Function}: A subset $\mathcal{A}%
_{a}\subset\mathcal{M}$ is said to be \emph{objective} if for every
$\mathbf{A}\in\mathcal{A}_{a}$ and every $\mathbf{Q}\in\mathcal{M}_{ort}^{+}$,
$\mathbf{QA}\in\mathcal{A}_{a}$. A scalar-valued function $U:\Omega
\times\mathcal{A}_{a}\rightarrow\mathbb{R}$ is said to be \emph{objective} if
its domain is objective and

$U(\mathbf{X,QA})=U(\mathbf{X,A})\quad\forall\mathbf{A}\in\mathcal{A}%
_{a},\quad\forall\mathbf{Q}\in\mathcal{M}_{ort}^{+}.\quad\quad(6.20)$

(D2) \emph{Isotropic Set and Isotropic Function} . A subset $\mathcal{A}%
_{a}\subset\mathcal{M}$ is said to be isotropic if for every $\mathbf{A}%
\in\mathcal{A}_{a}$ and every $\mathbf{Q}\in\mathcal{M}_{ort}^{+}$ also
$\mathbf{AQ}^{T}\in\mathcal{A}_{a}$. A scalar valued function $U:\Omega
\times\mathcal{A}_{a}\rightarrow\mathbb{R}$ is said to be \emph{isotropic} if
its domain is isotropic and

$U(\mathbf{X,AQ}^{T})=U(\mathbf{X,A})\quad\forall\mathbf{A}\in\mathcal{A}%
_{a},\quad\forall\mathbf{Q}\in\mathcal{M}_{ort}^{+}.\quad\quad(6.21)$%
\textquotedblright.

\smallskip

The set $\mathcal{M}$ appearing in \cite[Def.\ 6.1.2]{Gao:00}, as well as a
set $\mathcal{A}$, are defined on page 287 of \cite{Gao:00}:

\textquotedblleft We use the notation $\mathcal{M}(\Omega;\mathbb{R}^{m\times
n})$, or simply $\mathcal{M}$, to denote the space of all second-order tensor
functions with domain $\Omega$ in $\mathbb{R}^{n}$ and range in $\mathbb{R}%
^{m\times n}$. Let $\mathcal{A\subset M}$ be an admissible deformation
gradient space, defined by

$\mathcal{A}=\{\mathbf{A}\in\mathcal{M}\mid\operatorname*{rank}\mathbf{A}%
(\mathbf{X})=\min\{m,n\}\ \ \forall\mathbf{X}\in\Omega\}.\quad\quad
(6.16)$\textquotedblright.

\smallskip

For us Definition 6.1.2 is not sufficiently clear because $\Omega
\times\mathcal{A}_{a}$ does not seem to be an objective set. This is the
reason to use Definition 1 for our characterization of objective sets and
objective functions.

\smallskip

Observe that in the first version of \cite{GaoWu:11v2} one finds another
definition of an objective function which doesn't seem to be equivalent to
that in Definition 1 above. Indeed, in \cite{GaoWu:11v1} one can read:

\textquotedblleft By the fact that this potentially useful theory is based on
certain fundamental principles in physics, the nonconvex term $W(x)$ is
required to be an \emph{objective function}, i.e., there exists a
\emph{geometrically nonlinear mapping} $\Lambda:\mathcal{X}_{a}\rightarrow
\mathcal{V}\subset\mathbb{R}^{m}$ and a canonical function $V:\mathcal{V}%
\subset\mathbb{R}^{m}\rightarrow\mathbb{R}$ such that $W(x)=V(\Lambda(x))$
$\forall x\in\mathcal{X}_{a}$.\textquotedblright

Almost the same text can be found in \cite[p.\ 118]{GaoWu:11-mame}%
.\footnote{This paper is contained in the file
\url{http://www.confmame.ru/documents/pdf/p1MAME'IV.pdf}.}

\section{Characterizations of objective sets and objective functions}

In the sequel $\mathbb{R}^{m}$ $(m\geq1)$ is endowed with the usual inner
product, and the elements of $\mathbb{R}^{m}$ are considered as being
$m\times1$ matrices. Let set
\[
\mathcal{Q}_{m}:=\left\{  Q\in\mathbb{R}^{m\times m}\mid Q^{T}Q=I_{m},\ \det
Q=1\right\}  .
\]

\begin{lemma}
\label{lem1}Let $m\geq2$ and $u,v\in S^{m-1}:=\{x\in\mathbb{R}^{m}%
\mid\left\Vert x\right\Vert =1\}$. Then there exists $Q\in\mathcal{Q}_{m}$
such that $Qu=v$.
\end{lemma}

Proof. Consider first the case $m=2$. Then there exist $\alpha,\beta\in
\lbrack0,2\pi)$ such that $u=[\cos\alpha~\sin\alpha]^{T}$, $v=[\cos\beta
~\sin\beta]^{T}$. Taking
\[
Q_{2}=\left[
\begin{array}
[c]{cc}%
\cos(\beta-\alpha) & -\sin(\beta-\alpha)\\
\sin(\beta-\alpha) & \cos(\beta-\alpha)
\end{array}
\right]  ,
\]
clearly $Q_{2}\in\mathcal{Q}_{2}$ and $v=Q_{2}u$. If $m\geq3$, take
$X\subset\mathbb{R}^{m}$ a linear subspace containing $u,v$, then consider a
orthonormal basis $e_{1},e_{2},\ldots,e_{m}$ in $\mathbb{R}^{m}$ with
$e_{1},e_{2}\in X$. With respect to this basis take $Q$ of the form
\[
Q=\left[
\begin{array}
[c]{cc}%
Q_{2} & 0\\
0 & I_{m-2}%
\end{array}
\right]  ,
\]
where $Q_{2}$ is defined as above and $I_{k}$ is the identity matrix of order
$k$. It is clear that $Q\in\mathcal{Q}_{m}$ and $Qu=v$.

\begin{proposition}
\label{prop1}Let $m\geq2$, $A\subset\mathbb{R}^{m}$ a nonempty set, and
$f:A\rightarrow E$.

\emph{(i)} $A$ is objective if and only if there exists $\Gamma\subset
\mathbb{R}_{+}$ a nonempty set such that $A=\Gamma\cdot S^{m-1}:=\{tx\mid
t\in\Gamma$, $x\in S^{m-1}\}$.

\emph{(ii)} $f$ is objective if and only if there exist $\Gamma\subset
\mathbb{R}_{+}$ a nonempty set and a function $\varphi:\Gamma\rightarrow E$
such that $A=\Gamma\cdot S^{m-1}$ and $f(x)=\varphi(\left\Vert x\right\Vert )$
for every $x\in A$.

\emph{(iii)} Every nonempty subset $B$ of $\mathbb{R}$ is objective and any
function $g:B\rightarrow E$ is objective.
\end{proposition}

Proof. (i) Let first have $A=\Gamma\cdot S^{m-1}$ with $\Gamma\subset
\mathbb{R}_{+}$ nonempty. Consider $x\in A$ and $Q\in\mathcal{Q}_{m}$. Then
$x=\gamma u$ with $\gamma\in\Gamma$ and $u\in S^{m-1}$. It follows that
$Qx=\gamma Qu$. Since $Q^{T}Q=I_{m}$, we get $\left\Vert Qu\right\Vert
^{2}=(Qu)^{T}Qu=u^{T}Q^{T}Qu=u^{T}I_{m}u=u^{T}u=\left\Vert u\right\Vert
^{2}=1$, whence $Qu\in S^{m-1}$. Hence $Qx\in\Gamma\cdot S^{m-1}=A$. Assume
now that $A$ is objective, and set $\Gamma:=\{\left\Vert x\right\Vert \mid
x\in A\}\subset\mathbb{R}_{+}$. Clearly, $A\subset\Gamma\cdot S^{m-1}$.
Consider $x\in\Gamma\cdot S^{m-1}$, that is $x=\gamma u$ with $\gamma\in
\Gamma$ and $u\in S^{m-1}$. Since $\gamma\in\Gamma$, there exists $y\in A$
such that $\gamma=\left\Vert y\right\Vert $. If $\gamma=0$ then $0=y\in A$,
whence $x=\gamma\cdot0=0\in A$. Assume that $\gamma\neq0$. Then $v:=\gamma
^{-1}y\in S^{m-1}$. By Lemma \ref{lem1}, there exists $Q\in\mathcal{Q}_{m}$
such that $u=Qv$. It follows that $x=\gamma u=\gamma Qv=Qy$. Since $A$ is
objective, it follows that $x\in A$. Hence $A=\Gamma\cdot S^{m-1}$.

(ii) Assume first that there exist $\Gamma\subset\mathbb{R}_{+}$ a nonempty
set with $A=\Gamma\cdot S^{m-1}$ and a function $\varphi:\Gamma\rightarrow E$
such that $f(x)=\varphi(\left\Vert x\right\Vert )$ for every $x\in A$; hence
$A$ is objective by (i). Since $\left\Vert Qx\right\Vert =\left\Vert
x\right\Vert $ for every $Q\in\mathcal{Q}_{m}$, it is clear that $f$ is objective.

Conversely, assume that $f$ is objective. By the very definition, $A$ is
objective, and so, by (i), there exists $\Gamma\subset\mathbb{R}_{+}$ a
nonempty set such that $A=\Gamma\cdot S^{m-1}$. Fix some $u_{0}\in S^{m-1}$
and take $\varphi:\Gamma\rightarrow E$ defined by $\varphi(t):=f(tu_{0})$.
Consider $x\in A$; then $\gamma:=\left\Vert x\right\Vert \in\Gamma$. If
$\gamma=0$, then $x=0=\gamma u_{0}$, and so $f(x)=f(\gamma u_{0}%
)=\varphi(\gamma)=\varphi(\left\Vert x\right\Vert )$. Assume that $\gamma>0$.
Then $u:=\gamma^{-1}x\in S^{m-1}$. By Lemma \ref{lem1} there exists
$Q\in\mathcal{Q}_{m}$ such that $u=Qu_{0}$, whence $x=Q(\gamma u_{0})$. Since
$f$ is objective, it follows that $f(x)=f\left(  Q(\gamma u_{0})\right)
=f(\gamma u_{0})=\varphi(t)=\varphi(\left\Vert x\right\Vert )$.

(iii) Clearly, $\mathcal{Q}_{1}=\left\{  [1]\right\}  =\left\{  I_{1}\right\}
$. The conclusion follows.

\begin{corollary}
\label{cor2}Let $m\geq2$ and $H\in\mathbb{R}^{m\times m}$. Consider
$A\subset\mathbb{R}^{m}$ a nonempty set such that $A\neq\{0\}$ and
$f:A\rightarrow\mathbb{R}$ defined by $f(x):=x^{T}Hx$. Then $f$ is objective
if and only if there exist $\emptyset\neq\Gamma\subset\mathbb{R}_{+}$ and
$\alpha\in\mathbb{R}$ such that $A=\Gamma\cdot S^{m-1}$ and $\tfrac{1}%
{2}(H+H^{T})=\alpha I_{m}$.
\end{corollary}

Proof. Setting $H_{s}:=\tfrac{1}{2}(H+H^{T})$, we have that $f(x):=x^{T}%
H_{s}x$ for $x\in A$. So we may (and do) assume that $H=H_{s}$. The
sufficiency is obvious because, by Proposition \ref{prop1} (i), $A$ is
objective and $f(x)=\alpha\left\Vert x\right\Vert ^{2}$ for $x\in A$. Hence
$f$ is objective by Proposition \ref{prop1} (ii) (with $\varphi(t)=\alpha
t^{2}$).

Assume that $f$ is objective. By Proposition \ref{prop1} (ii), there exist
$\emptyset\neq\Gamma\subset\mathbb{R}_{+}$ and $\varphi:\Gamma\rightarrow
\mathbb{R}$ such that $A=\Gamma\cdot S^{m-1}$ and $f(x)=\varphi(\left\Vert
x\right\Vert )$ for $x\in A$. Fixing $u_{0}\in S^{m-1}$, for $t_{0}\in
\Gamma\setminus\{0\}$ we have that
\[
\varphi(t_{0})=f(t_{0}u_{0})=f(t_{0}u)=(t_{0}u)^{T}H(t_{0}u)=t_{0}^{2}%
u^{T}Hu=\alpha t_{0}^{2}\quad\forall u\in S^{m-1},
\]
where $\alpha:=u_{0}^{T}Hu_{0}$. It follows $u^{T}Hu=\alpha u^{T}I_{m}u$ for
every $u\in\mathbb{R}^{m}$, and so $H_{s}=H=\alpha I_{m}$. The proof is complete.

\section{Further remarks}

1) On the webpage \url{https://en.wikipedia.org/wiki/Radial\_function} one
finds the following definition: \textquotedblleft a radial function is a
function defined on a Euclidean space $\boldsymbol{R}^{n}$ whose value at each
point depends only on the distance between that point and the
origin\textquotedblright, continued by \textquotedblleft A function is radial
if and only if it is invariant under all rotations leaving the origin fixed.
That is, $f$ is radial if and only if $f\circ\rho=f$ for all $\rho\in SO(n)$,
the special orthogonal group in $n$ dimensions\textquotedblright. As reference
it is mentioned Stein~\& Weiss' book (\cite[p.~33]{SteWei:71}).

Clearly, $SO(n)$ is nothing else than $\mathcal{Q}_{n}$. This shows that
Proposition \ref{prop1} (ii) is a somewhat more precise description of radial functions.

\smallskip2) Latorre~\& Gao in \cite[p.~1765]{LatGao:16} and Jin~\& Gao
\cite[p.~169]{JinGao:16} say \textquotedblleft any convex quadratic function
is objective", while Gao, Ruan~\& Latorre in \cite[p.~NP9]{GaoRuaLat:16} and
\cite[p.~7]{GaoRuaLat:17} say \textquotedblleft Clearly, any convex quadratic
function $W(\varepsilon)$ is objective due to the Cholesky decomposition
$A=\Lambda^{\ast}\Lambda\succeq0$."

\smallskip Therefore, these authors claim the following:

\smallskip\textbf{Claim 1} ({Gao et al.}): \emph{Any convex quadratic
function is objective}.

\smallskip To understand what D.Y. Gao and his collaborators mean by a
quadratic function let us quote from \cite[pp.~NP54, NP55]{LiuGaoWan:16} and
\cite[p.~160]{LiuGaoWan:17}:

\textquotedblleft In many applications, the geometrical operator $\Lambda$ is
usually a vector-valued quadratic function:

$\Lambda(x)=(\varepsilon(x),\phi(x))=\{\tfrac{1}{2}\left\langle x,B_{k}%
x\right\rangle +\left\langle x,b_{k}\right\rangle -d_{k},\tfrac{1}%
{2}\left\langle x,C_{s}x\right\rangle +\left\langle x,c_{s}\right\rangle
-e_{s}\}$,

\noindent where $b_{k}\in\mathbb{R}^{n}$ and $B_{k}\in\mathbb{R}^{n\times n}$
is a given symmetric matrix for each $k\in\{1,2,\cdots,p\}$; $c_{s}%
\in\mathbb{R}^{n}$ for each $s\in\{1,2,\cdots,q\}$, $C_{s}\in\mathbb{R}%
^{n\times n}$ is a given positive definite matrix for each $s\in
\{1,2,\cdots,q\}$; $d\in\mathbb{R}^{p}$ and $e\in\mathbb{R}^{q}$."

Consequently, this text shows that D.Y. Gao accepts (as myself) that a
quadratic function on $\mathbb{R}^{n}$ is a function $q:\mathbb{R}%
^{n}\rightarrow\mathbb{R}$ defined by $q(x):=\tfrac{1}{2}\left\langle
x,Cx\right\rangle +\left\langle x,d\right\rangle +e$ with $C\in\mathbb{R}%
^{n\times n}$ symmetric, $d\in\mathbb{R}^{n}$ and $e\in\mathbb{R},$
$\left\langle \cdot,\cdot\right\rangle $ being the usual inner product on
$\mathbb{R}^{n}$.

Using Claim 1 we get

\smallskip

\textbf{Claim 2}: \emph{Any linear function on }$\mathbb{R}^{n}$\emph{ is
objective}.

\smallskip Of course, Claim 2 contradicts the following assertion from
\cite[p.~3]{Gao:18}: \textquotedblleft an objective function must be nonlinear".

\smallskip

Coming back to Claim 1, consider $n\geq2$ and the function $q_{d}%
:\mathbb{R}^{n}\rightarrow\mathbb{R}$ defined by $q_{d}(x):=\left\Vert
x\right\Vert ^{2}+\left\langle x,d\right\rangle $ for a fixed $d\in
\mathbb{R}^{n};$ clearly, $q_{d}$ is convex and quadratic (it is obtained
taking $C:=I_{n}$, $d\in\mathbb{R}^{n}$ and $e:=0$.\footnote{We could take
$C:=0_{n}$ (the null element of $\mathbb{R}^{n\times n}$).}

Assuming that Claim 1 is true, it follows that for all $Q\in\mathcal{Q}_{n}$
we have
\[
\left\Vert x\right\Vert ^{2}+\left\langle x,Q^{T}d\right\rangle =\left\Vert
Qx\right\Vert ^{2}+\left\langle Qx,d\right\rangle =q_{d}(Qx)=q_{d}%
(x)=\left\Vert x\right\Vert ^{2}+\left\langle x,d\right\rangle \quad\forall
x\in\mathbb{R}^{n},
\]
and so $\left\langle x,Q^{T}d\right\rangle =\left\langle x,d\right\rangle $
for $x\in\mathbb{R}^{n}$; hence $Q^{T}d=d$ for all $Q\in\mathcal{Q}_{n}$. Fix
$k\in\overline{2,n}$ and consider $Q_{k}\in\mathbb{R}^{n\times n}$ whose
elements are defined as follows: $q_{ij}:=0$ for $i,j\in\overline{1,n},$
$i\neq j;$ $q_{11}:=q_{kk}:=-1;$ $q_{ii}:=1$ for $i\in\overline{2,n}%
\setminus\{k\}$. Of course, $Q_{k}=Q_{k}^{T}$ and $\det Q_{k}=1;$ hence
$Q_{k}\in\mathcal{Q}_{n}$. From $Q_{k}^{T}d=d$ we get $2d_{1}=2d_{k}=0,$ and
so $d_{i}=0$ for every $i\in\overline{1,n}$. This shows that $d=0$. Therefore,
Claim 1 is false for $n\geq2$ and $d\neq0$.

In fact, using Corollary \ref{cor2}, we get:

\smallskip
\textbf{Fact 1}. \emph{A quadratic function }$q:\mathbb{R}
^{n}\rightarrow\mathbb{R}$ $(n\geq2)$ \emph{is ``objective" if and only
if there exists }$\alpha,\beta\in\mathbb{R}$ \emph{such
that }$q(x)=\alpha\left\Vert x\right\Vert ^{2}+\beta$\emph{ for all}
$x\in\mathbb{R}^{n}$.

\smallskip3) D.Y. Gao provides the following result in \cite{Gao:16}:

\smallskip\textquotedblleft\textbf{Lemma 1} A real-valued function $W(w)$ is
objective if and only if there exists a real-valued function $\Phi(E)$ such
that $W(w)=\Phi(w^{T}w)$",

\smallskip\noindent without proof or reference. Notice that the first version
of the present note (\url{http://arxiv.org/abs/1407.1241}) is mentioned in the
bibliography of the third version of \cite{Gao:16} as reference [68], but not
cited in connection with Lemma 1.

\smallskip4) Jin~\& Gao in \cite[p.~171]{JinGao:16} say:

\textquotedblleft According to the definition of the objectivity, a nonconvex
function $W:Y\rightarrow\mathbb{R}$ is objective if and only if there exists a
function $V:Y\times Y\rightarrow\mathbb{R}$ such that $W(y)=V(y^{T}y)$ [4, 17]"

Excepting ``[4, 17]", the same text can be found in
\cite[p.~205]{JinGao:17}. In \cite{JinGao:16}, ``[4]" and
``[17]" are our references \cite{Cia:13} and
\cite{Gao:16}, respectively.

It is clear that the reference ``[17]" [that is \cite{Gao:16}; see also
Remark 3) above] is not adequate for the representation
$W(y)=V(y^{T}y)$ of an ``objective function". In fact it
is not possible to speak about $V(y^{T}y)$ for a function $V:Y\times
Y\rightarrow\mathbb{R}$. The term $y^{T}y$ suggests that $y$ is an element of
$\mathbb{R}^{n}$ $(=\mathbb{R}^{n\times1})$ (for some positive integer $n$),
and so $y^{T}y\in\mathbb{R};$ surely, $y^{T}y$ cannot be an element of
$Y\times Y$.

\smallskip5) In \cite[p.~3]{Gao:18} Gao says:

\smallskip\textquotedblleft According to [10], $G(\mathbf{g})$ is an objective
function if and only if there exists an objective measure $\epsilon
=\mathbf{g}^{T}\mathbf{g}$ and a real-valued function $\Phi(\epsilon)$ such
that $G(\mathbf{g})=\Phi(\mathbf{g}^{T}\mathbf{g})$."

\smallskip The reference \textquotedblleft\lbrack10]" from the above text is
\cite{Cia:13}; moreover, the characterization of an \textquotedblleft
objective function" by using the representation $G(\mathbf{g})=\Phi
(\mathbf{g}^{T}\mathbf{g})$ is nothing else than \cite[Lem.~1]{Gao:16} [see
Remark 3) above].

\smallskip6) D. Y. Gao and his collaborators like to assert

\smallskip\textbf{Claim 3}. \emph{The objectivity is not an assumption, but an
axiom.}

\smallskip Claim 3 is attributed explicitly to Ciarlet \cite{Cia:13} in
\cite[p.~1765]{LatGao:16}: `It is emphasized in recent book by Ciarlet [3]
that \textquotedblleft the objectivity is not an assumption but an
axiom\textquotedblright;' see also: \cite[p.~6]{GaoRuaLat:14}, \cite[p.~3]%
{Gao:16}, \cite[p.~NP141]{CheGao:16}, \cite[p.~NP8]{GaoRuaLat:16},
\cite[p.~169]{JinGao:16}, \cite[p.~293]{CheGao:17}, \cite[p.~7]{GaoRuaLat:17}.
However, \cite[p.~NP197]{RuaGao:16} is an exception; here it is said:
\textquotedblleft It was emphasized in [10] that `the objectivity is not an
assumption, but an axiom'," where \textquotedblleft\lbrack10]" is
\textquotedblleft Ruan, N, and Gao, DY. Canonical duality approach for
nonlinear dynamical systems. IMA J Appl Math 2014; 79(2): 313--325". Notice
that \cite{Cia:13} is mentioned in \cite{RuaGao:16} as reference
\textquotedblleft\lbrack20]". Probably \textquotedblleft\lbrack10]" in
\textquotedblleft It was emphasized in [10] ..." is a misprint, having in view
that the string \textquotedblleft object" can be found only two times in
reference \textquotedblleft\lbrack10]" of \cite{RuaGao:16}; it appears in the
text \textquotedblleft The optimal objective function value" from pages 322
and 323 of \textquotedblleft\lbrack10]".

\smallskip Having in view the texts reproduced from \cite[p.~1765]{LatGao:16}
and \cite[p.~3]{Gao:18}, we hoped to find the strings \textquotedblleft
objectivity" and \textquotedblleft objective function" in the book
\cite{Cia:13}. We were very disappointed because the string \textquotedblleft
objectiv" can be found in \cite{Cia:13} only in the word \textquotedblleft
objective", and this one used with its usual meaning of \textquotedblleft aim"
or \textquotedblleft scope", at pages: 1, 76, 193, 312, 394, 412, 437, 572,
606, 653, 675, 682, 711, 763 (The objective of this), 50 (our main objective
in this), 61 (the first objective of this), 64, 446, 492 (Our objective is
to), 98, 514 (The objective is to), 142, 423 (Our objective then consists),
241 (the objective consists), 248 (The "ideal objective" in this respect; this
objective is unattainable), 249 (the above "ideal objective"), 479 (these
objectives will be achieved), 523 (Our basic objective is to), 540 (Our first
objective is to), 637 (Our objective in this), 657 (Our more modest objective
in this), 726 (two objectives are achieved).

\smallskip7) As seen in Remarks 4) and 5), the representation of
\textquotedblleft objective functions" as in \cite[Lem.~1]{Gao:16} is
attributed to \cite{Cia:13}, too. Something related to such a representation
can be found at page 695 of \cite{Cia:13}:

\smallskip\textquotedblleft Remarks (1) It can be shown that the \emph{axiom
of material frame-indifference} implies that, as a function of $\mathbf{F}%
\in\mathbb{M}_{+}^{3}$, the stored energy function is in fact a function of
$\mathbf{F}^{T}\mathbf{F}\in\mathbb{S}_{>}^{3}$, where $\mathbb{S}_{>}^{3}$
denotes the set of all symmetric, positive-definite, symmetric matrices of
order three. In other words, there exists a function $\widetilde{W}%
:\overline{\Omega}\times\mathbb{S}_{>}^{3}\rightarrow\mathbb{R}$ such that, at
each $x\in\overline{\Omega}$, $W(x,\mathbf{F})=\widetilde{W}(x,\mathbf{F}%
^{T}\mathbf{F})$ for all $\mathbf{F}\in\mathbb{M}_{+}^{3}$. It therefore
follows that $W(x,\nabla\psi(x))=\widetilde{W}(x,\nabla\psi(x)^{T}\nabla
\psi(x))$ at each $x\in\overline{\Omega}$, where $\nabla\psi(x)^{T}\nabla
\psi(x)\in\mathbb{S}^{3}$ is none other than the \emph{metric tensor} at $x$
associated with the deformation $\psi$ (Section 8.2), also called in
elasticity theory the \emph{Cauchy-Green strain tensor} at $x$."

\smallskip So, it is not possible to attribute to \cite{Cia:13} the
representation of a radial function from \cite[p.~171]{JinGao:16} or from
\cite[p.~3]{Gao:18}, at least because the above remark from \cite[p.~695]%
{Cia:13} refers to elements $\mathbf{F}\in\mathbb{M}_{+}^{3}$ and not to
elements $y$ belonging to an ``objective"  subset of $\mathbb{R}^n$ (as the domain
of an \textquotedblleft objective function").

\end{document}